\documentclass[11pt,amsfonts]{article}
\usepackage{graphicx}
\usepackage{latexsym}
\usepackage{amssymb}
\usepackage{amsmath}
\usepackage{enumerate}

\usepackage{layout}
\usepackage{eufrak}
\newtheorem{prop}{Proposition}
\newtheorem{lemma}{Lemma}

\newtheorem{corollary}{Corollary}
\newtheorem{theorem}{Theorem}
\newtheorem{remark}{Remark}

\def\real{{\mathord{{\rm I\kern-2.8pt R}}}}        
\def\inte{{\mathord{{\rm I\kern-2.8pt N}}}}

\def\sZZ{{\rm Z\kern-2.8ptem{}Z}}

\def\z{{\mathchoice
  {\sZZ}
  {\sZZ}
  {\rm Z\kern-0.30em{}Z}
  {\rm Z\kern-0.25em{}Z} }}
\def\sQQ{{\kern 0.27em \vrule height1.45ex width0.03em depth0em
          \kern-0.30em \rm Q}}
\def\qu{{\mathchoice
    {\sQQ}
    {\sQQ}
  {\kern 0.225em \vrule height1.05ex width0.025em depth0em \kern-0.25em \rm Q}
  {\kern 0.180em \vrule height0.78ex width0.020em depth0em \kern-0.20em \rm Q}
        }}
\def\sCC{{\kern 0.27em \vrule height1.45ex width0.03em depth0em
          \kern-0.30em \rm C}}
\def\complex{{\mathchoice
    {\sCC}
    {\sCC}
  {\kern 0.225em \vrule height1.05ex width0.025em depth0em \kern-0.25em \rm C}
  {\kern 0.180em \vrule height0.78ex width0.020em depth0em \kern-0.20em \rm C}
        }}

\newcommand{\ba}{\begin{array}}
\newcommand{\ea}{\end{array}}
\newcommand{\be}{\begin{equation}}
\newcommand{\ee}{\end{equation}}
\newcommand{\bea}{\begin{eqnarray}}
\newcommand{\eea}{\end{eqnarray}}
\newcommand{\beaa}{\begin{eqnarray*}}
\newcommand{\eeaa}{\end{eqnarray*}}

\def\z{\zeta}

\font\tenmath=msbm10 \font\sevenmath=msbm7 \font\fivemath=msbm5
\newfam\mathfam \textfont\mathfam=\tenmath
\scriptfont\mathfam=\sevenmath \scriptscriptfont\mathfam=\fivemath

\def \={{\buildrel {\rm (law)} \over =}}

\def\cB{{\cal B}}

\def\cE{{\cal E}}
\def\cF{{\cal F}}

\def\cH{{\cal H}}

\def\cP{{\cal P}}

\def\cS{{\cal S}}

\def\qed{ \hfill \vrule width.25cm height.25cm depth0cm\smallskip}

\newcommand{\basa}{\begin{assumption}}
\newcommand{\easa}{\end{assumption}}

\newcommand{\bas}{\begin{assum}}
\newcommand{\eas}{\end{assum}}

\newcommand{\ignore}[1]{}
\begin{document}

\renewcommand{\thefootnote}{\fnsymbol{footnote}}

\renewcommand{\thefootnote}{\fnsymbol{footnote}}

\title{Hitting times for the stochastic wave equation with fractional-colored noise}
\author{Jorge Clarke De la Cerda $^{1,2}$ \footnote{Partially supported by the MECESUP proyect UCO-0713 and the CONICYT-ECOS program C10E03.} $\qquad $
Ciprian A. Tudor $^{2,}$ \footnote{Associate member of the team Samm, Universit\'e de Paris 1 Panth\'eon-Sorbonne. Partially supported by the ANR grant "Masterie" BLAN 012103. }\vspace*{0.1in} \\
$^{1}$ CI$^{2}$MA, Departamento de Ingenier\'ia Matem\'atica, Universidad de Concepci\'on,\\
Casilla 160-C, Concepci\'on, Chile. \\
jclarke@udec.cl \vspace*{0.1in} \\
 $^{2}$ Laboratoire Paul Painlev\'e, Universit\'e de Lille 1\\
 F-59655 Villeneuve d'Ascq, France.\\
 \quad tudor@math.univ-lille1.fr\vspace*{0.1in}}

\maketitle

\begin{abstract}
We give sharp regularity results for the solution to the stochastic wave equation with linear fractional-colored noise.  We apply these results in order to establish upper and lower bound for the hitting probabilities of the solution in terms of the Hausdorff measure and of the  Newtonian capacity.
\end{abstract}



\vskip0.3cm	


\section{Introduction}
\hspace{1.1cm} The recent development of the stochastic calculus with respect to the fractional Brownian motion (fBm) naturally led to the study of stochastic partial differential equations (SPDEs) driven by this Gaussian process. The motivation comes from the wide area of applications of the fBm. We refer, among others, to \cite{GLT06}, \cite{maslowski-nualart03}, \cite{nualart-vuillermont06}, \cite{QS-tindel07} and \cite{TTV}.
The purpose  of our paper is to study the stochastic wave equation driven by fractional-colored Gaussian noise. Our work is situated somehow in the continuation of the line of research which concerns SPDEs driven by the fBm but in the same time it follows the l research initiated by Dalang in \cite{dalang99} which treats equations with white noise in time and  correlated in space.
More precisely, we consider   a system of $k$ stochastic wave equations
\begin{equation}
\label{system}
\frac{\partial ^{2} u_{i} }{\partial t^{2} }(t,x)= \Delta u_{i}(t,x) + \dot{W}_{i}(t,x), \hskip0.5cm t\in [0,T], x\in \mathbb{R} ^{d}
\end{equation}
with initial condition $u_{i}(t,x)=0$ and $\frac{\partial u _{i}}{\partial t} (0,x)= 0$ for every $x\in \mathbb{R} ^{d}$ and for every $i=1,...,k$.
The driving Gaussian process behaves as a fractional Brownian motion in time and has spatial covariance given by the Riesz kernel. More precisely
\begin{equation*}
\mathbf{E} (W_{i}(t,A) W_{j}(s,B))= \delta _{i,j}R_{H}(t,s) \int_{A} \int_{B} f(x-y) dxdy
\end{equation*}
for every $t,s \in [0,T]$ and $A,B$ Borel sets in $\mathbb{R} ^{d}$ where $f:\mathbb{R}^{d} \to \mathbb{R} _{+}$ is the Fourier transform  of a non-negative tempered measure $\mu$ on $\mathbb{R} ^{d}$ whose density with respect to the Lebesque measure is $\vert \xi \vert ^{-(d-\beta)}$, $0<\beta <d$. Above $\delta_{i,j}$ denotes  the Kronecker symbol.

The equation (\ref{system}) has been recently studied in \cite{BT}. It is has been proven that (\ref{system}) admits a unique mild solution if and only if $\beta < 2H+1$ which extends the result obtained in \cite{dalang99} in the case $H=\frac{1}{2}$. The purpose of this work is to analyze further the solution of (\ref{system}).  We will actually  give sharp results for the regularity of it, in time and in space, and we apply these regularity results to study the hitting probabilities for the solution $u$ to (\ref{system}). More precisely, given a Borel set $A \subset \mathbb{R} ^{k}$ we want to determine whether the process $(u(t,x), t\in [0,T], x\in \mathbb{R} ^{d})$ hits the set $A$ with positive probability. Recently, there has been several papers on hitting probabilities, and more generally speaking, on potential theory for systems of SPDEs. We refer, among others, to  \cite{DaNu}, \cite{DKN1}, \cite{DKN2}, \cite{DaSa2} or \cite{MuTr}. The study of hitting probabilities for stochastic partial differential equations with fractional noise in time is new. As far as we know, only the paper \cite{EF} treated  this problem. Actually, in this reference the authors give upper and lower bounds for the hitting times of  solution to a system of stochastic heat equations on the circle with fractional noise in time.

Our aim is to make a new step in this research direction. As we mentioned before we make a potential analysis of the solution to the stochastic wave equation with fractional-colored noise. That means, the noise behaves as the fractional Brownian motion with respect to the time variable and it is a "colored" non-white spatial covariance. In our work this spatial covariance will be described by the Riesz kernel.  It is know classical the fact that in order to obtain results on the hitting times of a stochastic process, a detailed analysis of the behavior of the increments of the process is needed. We address this question in our paper and we find the following: the solution $u(t,x), t\in [0,T], x\in \mathbb{R} ^{d}$ to (\ref{system}) is H\"older continuous of order $2H+1-\beta$, $\beta \in (2H-1, d\wedge 2H+1)$  in time as well as with respect to the  space variable. This generalizes the result obtained in \cite{DaSa2} and \cite{DaSa1} for the wave equation with white noise in time and Riesz covariance in space. Although the main lines of our work follows the approach  of \cite{DaSa1}, we stress that, as usually, the fractional cases involves more complex calculation and the techniques used in the standard white noise case need to be substantially  adapted. this is mainly due to the nature of the noise and to the structure of the Gaussian space associated to the noise. We will point out later in our paper, how the  fractional noise  involves more complexity  than  in e.g. \cite{DKN1} or \cite{DaSa1}. Moreover, the study of the solution to the wave equation is generally recognized to be more difficult that e.g. the solution to the heat equation, due to the appearance of the trigonometric functions and this is also the case in our work.

   We mention that there are  more or less general criteria to determine the hitting times for a stochastic process. Such criteria have been given in \cite{BLX}, \cite{DKN1}, \cite{DKN2} or \cite{DaSa1} among others. We will use the approach in \cite{BLX} because it concerns Gaussian processes and fits well with our context (note that the solution to (\ref{system}) is Gaussian).

Our paper is structured as follows. Section 2 contains some preliminaries, we briefly describe the basic properties of the Gaussian noise and its associated Hilbert space, we list the elements of the potential theory that we will use in our paper and we will recall some fact related to the solution to the stochastic wave equation with fractional-colored noise. In Section 3 we analyze the H\"older regularity of the solution with respect to its time and space variables. Section 4 is devoted to the study of the hitting probabilities for this solution, based on a criterium in \cite{BLX}.

\section{Preliminaries}
\hspace{1.1cm} This section is devoted to introduce the basic notion that we will need throughout the paper. We first introduce the canonical Hilbert space associated to the fractional-colored Gaussian noise. In the second part we present the basic elements related to the potential theory that intervene in the last section.

\subsection{The canonical Hilbert space}

\hspace{1.1cm} We denote by $C_0^{\infty}(\mathbb{R}^{d+1})$ the space of infinitely
differentiable functions on $\mathbb{R}^{d+1}$ with compact support, and
$\cS(\mathbb{R}^d)$ the Schwartz space of rapidly decreasing $C^{\infty}$
functions in $\mathbb{R}^d$. For $\varphi \in L^1(\mathbb{R}^d)$, we let $\cF
\varphi$ be the Fourier transform of $\varphi$:
$$\cF \varphi (\xi)=\int_{\mathbb{R}^d} e^{-i \xi \cdot x}\varphi (x)dx.$$

We begin by introducing the framework of \cite{dalang99}. 
Let $\mu$ be a non-negative tempered measure on $\mathbb{R}^d$, i.e. a
non-negative measure which satisfies:
\begin{equation} \label{mu-tempered}
\int_{\mathbb{R}^d}
\left(\frac{1}{1+ |\xi|^{2}} \right)^l \mu(d\xi)<\infty, \quad
 \mbox{for some} \ l >0.
\end{equation}

Since the integrand is non-increasing in $l$, we may assume that $l
\geq 1$ is an integer. Note that $1+|\xi|^2$ behaves as a constant
around $0$, and as $|\xi|^2$ at $\infty$, and hence
(\ref{mu-tempered}) is equivalent to:
\begin{equation} \label{mu-tempered-equiv}
\int_{|\xi| \leq 1}\mu(d\xi)<\infty, \quad
\mbox{and} \quad \int_{|\xi| \geq
1}\mu(d\xi)\frac{1}{|\xi|^{2l}}<\infty, \quad
 \mbox{for some integer} \ l \geq 1.
\end{equation}

Let $f: \mathbb{R}^d \to \mathbb{R}_{+}$ be the Fourier transform of $\mu$ in
$\cS'(\mathbb{R}^d)$, i.e.
$$\int_{\mathbb{R}^d}f(x)\varphi(x)dx=\int_{\mathbb{R}^d}\cF
\varphi(\xi)\mu(d\xi), \quad \forall \varphi \in \cS(\mathbb{R}^d).$$
Simple properties of the Fourier transform show that for any
$\varphi, \psi \in \cS(\mathbb{R}^d)$,
$$\int_{\mathbb{R}^d} \int_{\mathbb{R}^d} \varphi(x)f(x-y)\psi(y)dx dy=
\int_{\mathbb{R}^d}\cF \varphi(\xi) \overline{\cF \psi(\xi)}\mu(d\xi).$$

An approximation argument shows that the previous equality also
holds for indicator functions $\varphi=1_{A},\psi=1_{B}$, with $A,B
\in \cB_{b}(\mathbb{R}^d)$, where $\cB_b(\mathbb{R}^d)$ is the class of bounded
Borel sets of $\mathbb{R}^d$:
\begin{equation}
\label{Fourier-indicator} \int_A \int_B f(x-y)dx dy=\int_{\mathbb{R}^d}\cF
1_{A}(\xi) \overline{\cF 1_{B}(\xi)} \mu(d\xi).
\end{equation}

Now we introduce the \emph{fractional Brownian motion} (fBm) with Hurst index $H \in (0,1)$. This is a zero-mean Gaussian process $(B^{H}_{t}) _{t\in [0,T]}$ with covariance
\begin{equation*}
R_{H}(t,s) := \frac{1}{2} \left(t^{2H} + s^{2H} - |t-s|^{2H} \right), \hskip0.5cm t,s\in [0,T].
\end{equation*}

Let us denote by ${\cal{H}}$ the canonical Hilbert space associated with this Gaussian process. This canonical Hilbert space is defined as the closure of the linear space  generated by the indicator functions $1_{[0,t]}, t\in [0,T]$ with respect to the inner product
\begin{equation*}
\langle 1_{[0,t]}, 1_{[0,s]} \rangle _{{\cal{H}} } = R_{H}(t,s).
\end{equation*}

It is well known that for $H>1/2$ we have the expression
\begin{equation}\label{formula-RH}
R_{H}(t,s)= \alpha _{H} \int_{0} ^{t} \int_{0} ^{s} \vert u-v\vert ^{2H-2} dudv
\end{equation}
for every $s, t\in [0,T]$ with $\alpha _{H}:= H(2H-1)$. More generally, for $H>1/2$ and every $\psi, \phi \in \cH=\cH\left([0,T]\right)$ we have
\begin{equation}\label{H-inner-product}
\langle \psi , \phi \rangle _{\cal{H}} = \alpha _{H} \int_{0} ^{T} \int_{0} ^{T} \psi(u) \phi(v) \vert u-v\vert ^{2H-2} dudv
\end{equation}

As in \cite{balan-tudor08}, on a complete
probability space $(\Omega,\cF,P)$, we consider a zero-mean Gaussian process $W=\{W_t(A); t \geq 0, A \in \cB_{b}(\mathbb{R}^d)\}$ with covariance:
\begin{equation} \label{noise-covariance}
\mathbf{E}(W_t(A)W_s(B))= R_{H}(t,s) \int_{A} \int_{B} f(x-y)dx dy=: \langle 1_{[0,t] \times A}, 1_{[0,s] \times B} \rangle_{\cH \cP}.
\end{equation}

Let $\cE$ be the set of linear combinations of elementary functions
$1_{[0,t] \times A}$, $t \geq 0, A \in \cB_b(\mathbb{R}^d)$, and $\cH \cP$
be the Hilbert space defined as the closure of $\cE$ with respect to
the inner product $\langle \cdot , \cdot \rangle_{\cH \cP}$.
(Alternatively, $\cH \cP$ can be defined as the completion of
 $C_0^{\infty}(\mathbb{R}^{d+1})$, with respect to the inner product
$\langle \cdot, \cdot \rangle_{\cH \cP}$; see \cite{balan-tudor08}.)

The map $1_{[0,t] \times A} \mapsto W_t(A)$ is an isometry between
$\cE$ and the Gaussian space $H^{W}$ of $W$, which can be extended
to $\cH \cP$. We denote this extension by: $$\varphi \mapsto
W(\varphi)=\int_0^{\infty}\int_{\mathbb{R}^d} \varphi(t,x)W(dt,dx).$$

In the present work, we assume that $H>1/2$. Hence,
(\ref{formula-RH}) holds. From (\ref{Fourier-indicator}) and
(\ref{formula-RH}), it follows that for any $\varphi,\psi \in \cE$,
\begin{eqnarray*}
\langle \varphi, \psi \rangle_{\cH \cP}&=& \alpha_H
\int_{0}^{\infty} \int_0^{\infty}\int_{\mathbb{R}^d}\int_{\mathbb{R}^d}
\varphi(u,x)\psi(v,y)f(x-y)|u-v|^{2H-2} dx \
dy \ du \ dv  \\
&=& \alpha_H  \int_{0}^{\infty} \int_0^{\infty} \int_{\mathbb{R}^d} \cF
\varphi(u,\cdot)(\xi)\overline{\cF\psi(v,\cdot)(\xi)} |u-v|^{2H-2}
\mu(d\xi)\ du \ dv.
\end{eqnarray*}

Moreover, we can interchange the order of the integrals $dudv$ and
$\mu(d\xi)$, since for indicator functions $\varphi$ and $\psi$, the
integrand is a product of a function of $(u,v)$ and a function of
$\xi$. Hence, for $\varphi,\psi \in \cE$, we have:
\begin{equation}
\label{norm-HP-2} \langle \varphi, \psi \rangle_{\cH \cP}= \alpha_H
\int_{\mathbb{R}^d} \int_{0}^{\infty} \int_0^{\infty} \cF
\varphi(u,\cdot)(\xi)\overline{\cF\psi(v,\cdot)(\xi)} |u-v|^{2H-2}
du \ dv \ \mu(d\xi).
\end{equation}


The space $\cH \cP$ may contain distributions, but contains the
space $|\cH \cP|$ of measurable functions $\varphi: \mathbb{R}_{+} \times
\mathbb{R}^d \to \mathbb{R}$ such that
$$\|\varphi \|_{|\cH \cP|}^2:=\alpha_H
\int_{0}^{\infty} \int_0^{\infty}\int_{\mathbb{R}^d}\int_{\mathbb{R}^d}
|\varphi(u,x)||\varphi(v,y)|f(x-y)|u-v|^{2H-2} dx \ dy \ du \
dv<\infty.$$


\subsection{Elements of the potential theory}
\hspace{1.1cm} Our aim is to analyze the probability
 \begin{equation*}
 P \left( u(I) \cap A \right) \not= \emptyset
 \end{equation*}
 where $u$ is the solution to (\ref{system}), $I$ is a Borel set included in $[0,T] \times \mathbb{R}^{d}$ and $A$ is a Borel set in $\mathbb{R} ^{k}$.  Here $u(I)$ means the image of $I$ under the random map $(t,x)\to  u(t,x)$.

We will briefly present the notion of the potential theory that we will need in our paper. For all Borel sets $F \subset \mathbb{R} ^{d}$ we define ${\cal{P}}(F) $ to be the set of all probability measures with compact support included in $F$. For all $\mu \in {\cal{P}} (\mathbb{R} ^{d})$, let us denote by $I_{\beta} (\mu)$ the so-called $\beta$-energy of the measure $\mu $ defined by
\begin{equation}
\label{i}
I_{\beta}(\mu)= \int _{\mathbb{R}^{d}}\int _{\mathbb{R}^{d}}K_{\beta} (\Vert x-y\Vert )\mu (dx)\mu (dy)
\end{equation}
where
\begin{equation}\label{Riesz-Kernel}
K_{\beta}(r)=
\begin{cases}
r^{-\beta} & \text{if $\beta>0$;}\\
log\left( \frac{N_{0}}{r} \right) & \text{if $\beta=0$;}\\
1 & \text{if $\beta<0$.}
\end{cases}
\end{equation}
Here $N_{0}$ is a constant.

For all $\beta \in \mathbb{R}$ and $F\in {\cal{B}} (\mathbb{R} ^{d})$ we define the $\beta$-dimensional capacity of $F$ by
\begin{equation}\label{Riesz capacity}
\text{Cap} _{\beta } (F)= \left[ \inf _{\mu \in {\cal{P}}(F) } I_{\beta } (\mu) \right] ^{-1}
\end{equation}
with the convention $1/\infty:=0$. The $\beta$-dimensional Hausdorff measure of the set $F\in {\cal{B}} (\mathbb{R} ^{d})$ is given by
\begin{equation}\label{Hausdorff measure}
{\cal{H}} _{\beta } (F)= \lim _{\varepsilon \to 0^{+}} \inf \left[ \sum_{i=1}^{\infty} (2r_{i}) ^{\beta} ;\; F \subset \bigcup _{i=1}^{\infty} B(x_{i}, r_{i}),\; \sup_{i\geq 1} r_{i} \leq \varepsilon \right]
\end{equation}
where $B(x,r)$ denotes the Euclidean ball of radius $r>0$ centered at $x\in \mathbb{R} ^{d}$. When $\beta <0$, the $\beta$-dimensional Hausdorff measure of $F$ is infinite by definition.

\subsection{The  stochastic wave equation with linear fractional-colored noise}

 Consider the linear stochastic wave equation driven by an infinite-dimensional fractional Brownian motion $ W $ with Hurst parameter $H\in (0,1)$. That is
\begin{eqnarray}
\label{wave}
 \frac{\partial^2 u}{\partial t^2}(t,x)&=& \Delta u
(t,x) +\dot W(t,x), \quad t>0, x \in \mathbb{R} ^{d} \\
\nonumber
u(0, x)&=& 0, \quad x \in \mathbb{R} ^{d} \\
\nonumber \frac{\partial u}{\partial t}(0,x) &=& 0, \quad x \in
\mathbb{R} ^{d}.
\end{eqnarray}
Here $\Delta$ is the Laplacian on $\mathbb{R} ^{d}$ and $W=\{W_t(A); t \geq 0, A \in \cB_{b}(\mathbb{R} ^{d})\}$ is a centered Gaussian field with covariance
$$\mathbf{E}(W_t(A)W_s(B))=R_{H}(t,s) \int_{A} \int_{B} f(x-y)dx dy,$$
where $R_{H}$ is the covariance of the fractional Brownian motion and $f$ is the Riesz kernel.

Let $G_1$ be the fundamental solution of $u_{tt}-\Delta u=0$. It is
known that $G_1(t, \cdot)$ is a distribution in ${\cal{S}'}(\mathbb{R}^d)$ with
rapid decrease, and
\begin{equation}
\label{Fourier-G-wave} \cF
G_1(t,\cdot)(\xi)=\frac{\sin(t|\xi|)}{|\xi|},
\end{equation}
for any $\xi \in \mathbb{R}^d,t>0,d \geq 1$ (see e.g. \cite{treves75}). In particular,
\begin{eqnarray*}
G_1(t,x)&=&\frac{1}{2}1_{\{|x|<t\}}, \quad \mbox{if} \ d=1 \\
G_1(t,x)&=&\frac{1}{2 \pi}\frac{1}{\sqrt{t^2-|x|^2}}1_{\{|x|<t\}},
\quad \mbox{if} \ d=2 \\
G_1(t,x)&=&c_{d}\frac{1}{t}\sigma_t, \quad \mbox{if} \ d=3,
\end{eqnarray*}
where $\sigma_t$ denotes the surface measure on the 3-dimensional
sphere of radius $t$.

The solution of (\ref{wave}) is a square-integrable process
$u=\{u(t,x); t \geq 0, x \in \mathbb{R}^d\}$ defined by:
\begin{equation} \label{def-sol-wave} u(t,x)=\int_{0}^{t}
\int_{\mathbb{R}^d}G_1(t-s,x-y)W(ds,dy).
\end{equation}

By definition, $u(t,x)$ exists if and only if the stochastic integral above is well-defined, i.e. $g_{tx}:=G_1(t-\cdot,x-\cdot)
\in \cH \cP$. In this case, $\mathbf{E}|u(t,x)|^2 = \|g_{tx}\|_{\cH \cP}^2$.

The following result has been proved in \cite{BT}.

\begin{theorem}
The stochastic wave equation (\ref{wave}) admits an unique mild solution $(u(t,x) )_{t\in [0,T] , x\in \mathbb{R} ^{d} } $ if and only if
\begin{equation}
\label{1}
\int_{\mathbb{R} ^{d}} \left( \frac{1}{1+ |\xi|^{2}} \right) ^{H+ \frac{1}{2}} \mu(d\xi) < \infty.
\end{equation}

\end{theorem}

\begin{remark}
Note that (\ref{1}) is equivalent to:
\begin{equation}
\label{mu-tempered-eq} \int_{|\xi| \leq 1}\mu(d\xi)<\infty, \quad
\mbox{and} \quad \int_{|\xi| \geq 1}\mu(d\xi)\frac{1}{|\xi|^{2H+1}}<\infty.
  \end{equation}

\end{remark}

As  mentioned in the Introduction, we will consider throughout the paper that the spatial covariance of the noise $W$ is given by the Riezs kernel. That means the measure $\mu$ is
\begin{equation*}
d\mu (\xi)= \vert \xi \vert ^{-d+\beta } d\xi \mbox{ with } \beta \in (0,d).
\end{equation*}
In this case the kernel $f$ is given by
\begin{equation*}
f(\xi)= \vert \xi \vert ^{-\beta} \mbox{ with } \beta \in (0,d).
\end{equation*}
Note that in the case of the Riesz kernel, condition (\ref{1}) is equivalent to
\begin{equation}
\label{b1}
\beta \in (0, d\wedge (2H+1)).
\end{equation}

\begin{remark} Since $H>\frac{1}{2}$ and so $2H+1 \in (2,3)$, for dimension $d=1,2$  we have $\beta \in (0,d)$ while for $d\geq 3$ we have $\beta \in (0, 2H+1)$.
\end{remark}

\section{Regularity of the solution}

\subsection{Time regularity}
\hspace*{1.1cm} In this part we will focus our attention on the behavior of the increments of  the solution $u(t,x)$ with respect to the variable $t$. We will give upper and lower bounds for the $L^{2}$-norm of this increment. Usually, obtaining upper bounds is recognized to be easier than obtaining lower bounds, this is also the case in our work. Actually, to get the sharpness of the regularity of $u$ with respect to the time variable, we need to impose a stronger assumption than (\ref{b1}) on the parameters $\beta$ and $H$ (condition (\ref{b2}) below). This is due to the characteristics of the scalar product in the ${\cal{HP}}$.

We will start with the following useful lemma that gives an explicit expression for the ${\cal{H}}$ norm of the cosinus and sinus functions. These norms will widely appear further in our computations.

\begin{lemma}\label{sincos}
Let $f(x)= \cos (x) $ and $g(x)= \sin x$ for $x\in \mathbb{R}$. Then  for every $a,b\in \mathbb{R}$, $a<b$
\begin{equation*}
\Vert f 1_{(a,b) } \Vert ^{2} _{{\cal{H}}}=   \alpha _{H} \int_{0} ^{b-a} dv \cos (v) v^{2H-2} (b-a-v) +\alpha _{H} \cos (a+b)  \int_{0} ^{b-a} dv v^{2H-2}\sin (b-a -v).
\end{equation*}
and
\begin{equation*}
\Vert g 1_{(a,b) } \Vert ^{2} _{{\cal{H}}}=\alpha _{H} \int_{0} ^{b-a} dv \cos (v) v^{2H-2} (b-a-v) -\alpha _{H} \cos (a+b)  \int_{0} ^{b-a} dv v^{2H-2}\sin (b-a -v).
\end{equation*}
\end{lemma}
{\bf Proof: } Using the expression of the scalar product in the Hilbert space ${\cal{H}}$ and the trigonometric identities
\begin{eqnarray*}
\cos (u \pm v) &=& \cos u \cos v \mp \sin u \sin v  \quad \mbox{and} \\
sin(x) - sin(y) &=& 2cos(\frac{x+y}{2})sin(\frac{x-y}{2})
\end{eqnarray*}
we can write
\begin{eqnarray*}
\Vert f 1_{(a,b) } \Vert ^{2} _{{\cal{H}}}+ \Vert g 1_{(a,b) } \Vert ^{2} _{{\cal{H}}}
&=& \alpha _{H} \int_{a}^{b} \int_{a} ^{b} \vert u-v\vert ^{2H-2} \left( \cos u \cos v + \sin u \sin v\right) dudv\\
&=&  \alpha _{H} \int_{a}^{b}du \int_{a} ^{b} dv\vert u-v\vert ^{2H-2} \cos (u-v) \\
&=& 2 \alpha _{H} \int_{a}^{b} du  \int_{0} ^{u-a}dv  \cos (v) v^{2H-2} \\
&=& 2 \alpha _{H} \int_{0} ^{b-a} dv \cos (v) v^{2H-2} (b-a-v)
\end{eqnarray*}
where we made the change of variables $\tilde{v}=u-v$ in the integral $dv$ above and we computed the integral $du$. Similarly
\begin{eqnarray*}
\Vert f 1_{(a,b) } \Vert ^{2} _{{\cal{H}}}- \Vert g 1_{(a,b) } \Vert ^{2} _{{\cal{H}}}
&=& \alpha _{H} \int_{a}^{b} du\int_{a} ^{b} dv \vert u-v\vert ^{2H-2} \left( \cos u \cos v - \sin u \sin v\right) \\
&=&\alpha _{H} \int_{a}^{b} \int_{a} ^{b} \vert u-v\vert ^{2H-2}\cos (u+v)dudv
\end{eqnarray*}
and by the change of variable $\tilde{v}= u-v$ in the integral $dv$,
 \begin{eqnarray*}
\Vert f 1_{(a,b) } \Vert ^{2} _{{\cal{H}}}- \Vert g 1_{(a,b) } \Vert ^{2} _{{\cal{H}}}
&=&  2 \alpha _{H}  \int_{a}^{b} du \int_{0} ^{u-a} dv \cos (2u-v) v^{2H-2} \\
&=& 2 \alpha _{H}  \int_{0} ^{b-a} dv v^{2H-2} \int_{v+a}^{b} du\cos (2u -v) \\
&=&  \alpha _{H}  \int_{0} ^{b-a} dv v^{2H-2} \left( \sin (2b-v) - \sin (2a+ v) \right) \\
&=&2 \alpha _{H} \cos (a+b)  \int_{0} ^{b-a} dv v^{2H-2}\sin (b-a -v).
\end{eqnarray*}

\qed

\begin{remark}\label{rem3}
As a consequence of the Lemma \ref{sincos} we deduce the following
\begin{description}
\item{$i. $} For any $x>0$ the quantity $\int_{0} ^{x} v^{2H-2} \cos (v) (x-v)dv$ is positive (it is the sum of two norms).
\item{$ii. $} For every $a,b\in \mathbb{R}$, $a<b$
\begin{equation*}
\Vert f 1_{(a,b) } \Vert ^{2} _{{\cal{H}}} \leq 2 \alpha _{H} \int_{0} ^{b-a} dv \cos (v) v^{2H-2} (b-a-v)
\end{equation*}
\item{$iii. $}   For every $a,b\in \mathbb{R}$, $a<b$ \begin{equation*}
\Vert f 1_{(a,b) } \Vert ^{2} _{{\cal{H}}} \geq 2 \alpha _{H} \cos (a+b)  \int_{0} ^{b-a} dv v^{2H-2}\sin (b-a -v).
\end{equation*}
\end{description}

\end{remark}

Later, we use also the following lemma.
\begin{lemma}
\label{10m-1}
For every $a, b\in \mathbb{R}$ with $a<b$,
\begin{equation*}
\int_{a} ^{b} \int_{a}^{b} dudv \sin (u-v) \vert u-v \vert ^{2H-2}= 0.
\end{equation*}
\end{lemma}
{\bf Proof: } This follows from the trivial equality
$$\int_{a} ^{b} \int_{a}^{b} \sin(u)\cos (v) \vert u-v\vert ^{2H-2}dudv =\int_{a} ^{b} \int_{a}^{b} \sin(v)\cos (u) \vert u-v\vert ^{2H-2}dudv.$$
\qed

\vskip0.2cm

Concretely, we will  prove the following result concerning the time regularity of the solution to (\ref{wave}). We mention that, in the rest of our paper, $c,C...$ will denote generic positive constants that may change from line to line.

\begin{prop}\label{ptime}
Assume that
\begin{equation}\label{b2}
\beta \in (2H-1, d\wedge (2H+1) ).
\end{equation}
Let $t_{0},M>0$ and fix $x\in [-M,M]^{d}$. Then there exists a positive constants $c_{1}, c_{2}$ such that for every $s,t\in [t_{0}, T]$
\begin{equation*}
c_{1} \vert t-s\vert ^{2H+1-\beta} \leq \mathbf{E} \left| u(t,x)-u(s,x) \right| ^{2} \leq c_{2} \vert t-s\vert  ^{2H+1-\beta}.
\end{equation*}
\end{prop}
{\bf Proof: } Let $h>0$ and let us estimate the $L^{2}(\Omega)$-norm of the increment $u(t+h, x)-u(t,x)$.
Splitting the interval $[0,t+h]$ into the
intervals $[0,t]$ and $[t,t+h]$, and using the inequality $|a+b|^2
\leq 2(a^2+b^2)$, we obtain:
\begin{eqnarray}
\mathbf{E}|u(t+h,x)-u(t,x)|^2 & \leq & 2 \{\|(g_{t+h,x}-g_{t,x})1_{[0,t]} \|_{\cH
\cP}^2 +
 \|g_{t+h,x}1_{[t,t+h]} \|_{\cH \cP}^2\} \nonumber \\
 &=:& 2[E_{1,t}(h)+E_2(h)].\label{n10}
 \end{eqnarray}
The first summand can be handled in the following way.
 \begin{eqnarray*}
E_{1,t}(h) &= & \alpha_H \int_{\mathbb{R}^d} \mu(d\xi )\int_0^t \int_0^t dv dv
|u-v|^{2H-2}\cF (g_{t+h,x}-g_{tx})(u,\cdot)(\xi) \\
& &\times  \overline{\cF
(g_{t+h,x}-g_{tx})(v,\cdot)(\xi)}\\
&=& \alpha_H \int_{\mathbb{R}^d} \mu(d\xi)\int_0^t \int_0^t du dv |u-v|^{2H-2}
[\cF G_1(u+h,\cdot)(\xi)-\cF G_1(u,\cdot)(\xi)] \\
& &\times  \overline{\cF G_1(v+h,\cdot)(\xi)-\cF G_1(v,\cdot)(\xi)} \\
&=& \alpha_H \int_0^t \int_0^t du dv |u-v|^{2H-2} I_{h},
\end{eqnarray*}
where
\begin{eqnarray*}
	I_{h} &=& \int_{\mathbb{R} ^{d}} \mu(d\xi ) [\cF G_1(u+h,\cdot)(\xi)-\cF G_1(u,\cdot)(\xi)][\overline{\cF G_1(v+h,\cdot)(\xi)-\cF G_1(v,\cdot)(\xi)}]\\
	&=& \int_{\mathbb{R} ^{d}} \mu(d\xi ) \frac{(\sin ((u+h)|\xi|) -\sin(u|\xi|))}{|\xi|}\frac{(\sin ((v+h)|\xi|) -\sin(v|\xi|))}{|\xi|}.
\end{eqnarray*}
\\
Using the last trigonometric identity presented before we obtain
\\
  \begin{eqnarray*}
  E_{1,t}(h)&=&\alpha _{H} \int_{0} ^{t} \int_{0} ^{t} dudv \vert u-v\vert ^{2H-2} \int_{\mathbb{R} ^{d}}
  \mu(d\xi )\frac{ \sin (\frac{h\vert \xi \vert}{2} ) ^{2}}{\vert \xi \vert ^{2}}
  \cos (\frac{ (2u+h) \vert \xi \vert}{2})\cos (\frac{ (2v+h) \vert \xi \vert}{2})\\
  &=& c \alpha _{H} \int_{0} ^{t} \int_{0} ^{t} dudv \vert u-v\vert ^{2H-2} \int_{\mathbb{R} ^{d}} \frac{d \xi }{\vert \xi \vert ^{d-\beta + 2}}\sin (h\vert \xi \vert ) ^{2} \cos ((2u+h) \vert \xi \vert )\cos ((2v+h) \vert \xi \vert ),
  \end{eqnarray*}
  and by making the change of variables $\tilde{u}= (2u+h) \vert \xi \vert, \tilde{v}= (2v+h) \vert \xi \vert $,
  \begin{eqnarray}
  E_{1,t}(h)&=& c \cdot \alpha _{H}\int_{\mathbb{R} ^{d}}
  \frac{d\xi}{\vert \xi \vert ^{d-\beta +2H+2 }}  \sin (h\vert \xi \vert ) ^{2} \int_{h\vert \xi \vert } ^{(2t+h) \vert \xi \vert }  \int_{h\vert \xi \vert } ^{(2t+h) \vert \xi \vert } dudv \vert u-v\vert ^{2H-2} \cos u \cos v\nonumber \\
  &=& c \int_{\mathbb{R} ^{d}}
  \frac{d\xi}{\vert \xi \vert ^{d-\beta +2H+2 }}  \sin (h\vert \xi \vert ) ^{2} \Vert \cos ( \cdot ) 1_{( h\vert \xi \vert , (2t+h) \vert \xi \vert )}(\cdot ) \Vert ^{2} _{\cal{H}},
   \label{rel1}
  \end{eqnarray}
and using Lemma \ref{sincos},
\begin{eqnarray}
  E_{1,t}(h)&=& c \int_{\mathbb{R} ^{d}}
  \frac{d\xi}{\vert \xi \vert ^{d-\beta +2H+2 }}  \sin (h\vert \xi \vert ) ^{2}  \times \left[\int_{0} ^{2t\vert \xi \vert } \cos (v) v^{2H-2} (2t\vert \xi \vert -v) dv \right. \nonumber \\
   && \left. + \cos (2t\vert \xi \vert + 2h\vert \xi \vert ) \int_{0} ^{2t\vert \xi \vert } v^{2H-2} ( \sin (2t\vert \xi \vert -v)) \right] \nonumber \\
  &=& c \int_{\mathbb{R} ^{d}}
  \frac{d\xi}{\vert \xi \vert ^{d-\beta +2H+2 }}  \sin (h\vert \xi \vert ) ^{2} \times \left[ 2t\vert \xi \vert \int_{0} ^{2t\xi \vert } \cos (v) v^{2H-2} dv \right. \nonumber \\ && \left. -\sin (2t\vert \xi \vert ) (2t\vert \xi \vert ) ^{2H-1} + (2H-1) \int_{0} ^{2t\vert \xi \vert } \sin (v) v^{2H-2}dv \right. \nonumber \\
  && \left. + \cos (2t\vert \xi \vert + 2h\vert \xi \vert ) \int_{0} ^{2t\vert \xi \vert } v^{2H-2} ( \sin (2t\vert \xi \vert -v)) \right] \label{n11}
  \end{eqnarray}
  where we use integration by parts. By Remark \ref{rem3}, point ii. we have the upper bound
  \begin{eqnarray*}
  E_{1,t}(h)&\leq & c \cdot \alpha _{H}\int_{\mathbb{R} ^{d}}
  \frac{d\xi}{\vert \xi \vert ^{d-\beta +2H+2 }}  \sin (h\vert \xi \vert ) ^{2}  \\
  &&\times \left[ 2t\vert \xi \vert \int_{0} ^{2t\xi \vert } \cos (v) v^{2H-2} dv -\sin (2t\vert \xi \vert ) (2t\vert \xi \vert ) ^{2H-1}\right. \\
   &&\left. + (2H-1) \int_{0} ^{2t\vert \xi \vert } \sin (v) v^{2H-2}dv \right].
   \end{eqnarray*}
   We will treat separately the three summands above. Concerning the first one,
  \begin{eqnarray*}
 && \int_{\mathbb{R} ^{d}}
  \frac{d\xi}{\vert \xi \vert ^{d-\beta +2H+2 }}  \sin (h\vert \xi \vert ) ^{2}  2t\vert \xi \vert \int_{0} ^{2t \vert \xi \vert } \cos (v) v^{2H-2} dv \\
  &=& c_{t, H} h^{2H+1-\beta}\left|  \int_{\mathbb{R} ^{d}}
  \frac{d\xi}{\vert \xi \vert ^{d-\beta +2H+1 }}  \sin (\vert \xi \vert ) ^{2}\int_{0} ^{ \frac{2t\vert \xi \vert }{h} } \cos (v) v^{2H-2} dv \right|\\
  &\leq &  c_{t, H} h^{2H+1-\beta}  \int_{\mathbb{R} ^{d}}
  \frac{d\xi}{\vert \xi \vert ^{d-\beta +2H+1 }}  \sin (\vert \xi \vert ) ^{2}\left|\int_{0} ^{ \frac{2t\vert \xi \vert }{h} } \cos (v) v^{2H-2} dv \right|\\
  &\leq &  c_{t, H} h^{2H+1-\beta}
  \end{eqnarray*}
using condition (\ref{b1})  and the fact that the integral $\int_{0} ^{\infty} \cos (v) v^{2H-2} dv$ is convergent (this implies that the function $x\in [0, \infty) \to \int_{0}^{x} \cos (v) v^{2H-2} dv $ admits a limit at infinity and it is therefore bounded). On the other hand
\begin{eqnarray*}
&&\int_{\mathbb{R} ^{d}}
  \frac{d\xi}{\vert \xi \vert ^{d-\beta +2H+2 }}  \sin (h\vert \xi \vert ) ^{2} \sin (2t\vert \xi \vert ) (2t\vert \xi \vert ) ^{2H-1} = c_{t} h^{3-\beta } \int_{\mathbb{R} ^{d}}
  \frac{d\xi}{\vert \xi \vert ^{d-\beta +3 }} \sin (\vert \xi \vert)^{2} \sin \left( \frac{2t\vert \xi \vert }{h}\right)\\
  && = c_{t} h^{3-\beta } \int_{\vert \xi \vert \leq 1}
  \frac{d\xi}{\vert \xi \vert ^{d-\beta +3 }} \sin (\vert \xi \vert)^{2} \sin \left( \frac{2t\vert \xi \vert }{h}\right)+ c_{t} h^{3-\beta } \int_{\vert \xi \vert >1}
  \frac{d\xi}{\vert \xi \vert ^{d-\beta +3 }} \sin (\vert \xi \vert)^{2} \sin \left( \frac{2t\vert \xi \vert }{h}\right).
  \end{eqnarray*}
  The second part over the region  $\vert \xi \vert \geq 1$ is bounded by $ch^{3-\beta}$ simply by majorizing sinus by one. The second integral has a singularity for $\vert \xi \vert $ close to zero. Using that $sin(x) \leq x$ for all $x \geq 0$, we will bound it above by
  \begin{eqnarray*}
 && h^{3-\beta } \int_{\vert \xi \vert \leq 1}
  \frac{d\xi}{\vert \xi \vert ^{d-\beta +3 }} \sin (\vert \xi \vert)^{2} \sin \left( \frac{2t\vert \xi \vert }{h}\right)\\
  &&\leq c_{t} h^{3-\beta} \int_{\vert \xi \vert \leq 1}
  \frac{d\xi}{\vert \xi \vert ^{d-\beta +3 }} \vert \xi \vert ^{2} \left| \sin \left( \frac{2t\vert \xi \vert }{h}\right)\right| ^{2-2H} \left| \sin \left( \frac{2t\vert \xi \vert }{h}\right)\right| ^{2H-1}\\
  &&\leq c_{t}h^{2H+1-\beta} \int_{\vert \xi \vert \leq 1}
  \frac{d\xi}{\vert \xi \vert ^{d-\beta +2H-1 }}
\end{eqnarray*}
where we bounded $\left| \sin \left( \frac{2t\vert \xi \vert }{h}\right)\right| ^{2-2H}$ by $c_{t} (\vert \xi \vert h^{-1} )^{2-2H}$ and $ \left| \sin \left( \frac{2t\vert \xi \vert }{h}\right)\right| ^{2H-1}$ by 1. The last integral is finite since $\beta >2H-1$ (assumption (\ref{b2})).  Finally
\begin{eqnarray}
&&\int_{\mathbb{R} ^{d}}
  \frac{d\xi}{\vert \xi \vert ^{d-\beta +2H+2 }}  \sin (h\vert \xi \vert ) ^{2}\int_{0} ^{2t\vert \xi \vert } \sin (v) v^{2H-2}dv \nonumber \\
&&= h^{2H+2-\beta} \int _{\mathbb{R} ^{d}}
  \frac{d\xi}{\vert \xi \vert ^{d-\beta +2H+2 }}\sin (\vert \xi \vert ) ^{2}\int_{0} ^{\frac{2t\vert \xi \vert}{h }} \sin (v) v^{2H-2}dv\nonumber\\
&&= h^{2H+2-\beta} \int _{\vert \xi \vert \leq 1}
  \frac{d\xi}{\vert \xi \vert ^{d-\beta +2H+2 }}\sin (\vert \xi \vert ) ^{2}\int_{0} ^{\frac{2t\vert \xi \vert}{h }} \sin (v) v^{2H-2}dv\nonumber\\
&&+  h^{2H+2-\beta} \int _{\vert \xi \vert \geq 1}
  \frac{d\xi}{\vert \xi \vert ^{d-\beta +2H+2 }}\sin (\vert \xi \vert ) ^{2}\int_{0} ^{\frac{2t\vert \xi \vert}{h }} \sin (v) v^{2H-2}dv\nonumber\\
&&\leq h^{2H+2-\beta} \int _{\vert \xi \vert \leq 1}
  \frac{d\xi}{\vert \xi \vert ^{d-\beta +2H+2 }}\vert \xi \vert  ^{2}\int_{0} ^{\frac{2t \vert \xi \vert}{h }} \vert \sin v \vert  v^{2H-2}dv\nonumber\\
&&+  h^{2H+2-\beta} \int _{\vert \xi \vert \geq 1}
  \frac{d\xi}{\vert \xi \vert ^{d-\beta +2H+2 }}\int_{0} ^{\frac{2t\vert \xi \vert}{h }} \sin (v) v^{2H-2}dv.\nonumber\\
\end{eqnarray}
Using again the fact that $\int_{0} ^{\infty} \sin (v) v^{2H-2} dv$ is convergent is easy to see that the integral over the region $\vert \xi \vert \geq 1$ is bounded by $c_{t}h^{2H+2-\beta}$. For the intgral over $\vert \xi \vert \leq 1$ we make the change of variables $\tilde{v}=\frac{vh}{\xi}$ and we get
\begin{eqnarray*}
& & h^{3-\beta} \int _{\vert \xi \vert \leq 1}
  \frac{d\xi}{\vert \xi \vert ^{d-\beta +1}}\int_{0} ^{2t} \vert \sin \left(\frac{v\vert \xi \vert}{h} \right) \vert \; v^{2H-2}dv\\
& & = h^{3-\beta} \int_{\vert \xi \vert \leq 1}
  \frac{d\xi}{\vert \xi \vert ^{d-\beta +1 }} \int_{0} ^{2t} \left| \sin \left( \frac{v \vert \xi \vert }{h}\right)\right| ^{2-2H} \left| \sin \left( \frac{v \vert \xi \vert }{h}\right)\right| ^{2H-1} \; v^{2H-2}dv \\
  & & \leq c_{t} h^{2h+1-\beta} \int_{\vert \xi \vert \leq 1}
  \frac{d\xi}{\vert \xi \vert ^{d-\beta +2H -1 }},
\end{eqnarray*}
where we have made the same considerations as for the second summand in the decomposition of $E_{1,t}(h)$. In this way, we obtained the upper bound for the summand $E_{1,t}(h)$ in (\ref{n10})
  \begin{equation}
  \label{rel3}
  E_{1,t}(h) \leq C h^{2H+1-\beta}.
  \end{equation}

We study now the term $E_{2}(h)$ in (\ref{n10}) (its notation $E_{2}(h)$ instead of $E_{2,t}(h)$ is due to the fact that it does not depend on $t$, see below).  Using successively the change of variables $\tilde{u}= \frac{u}{h}, \tilde{v}=\frac{v}{h}$ in the integral $dudv$ and $\tilde{\xi}= h\xi$ in the integral $d\xi$,  the summand $E_{2} (h)$ can be written as
\begin{eqnarray*}
E_{2}(h)&=& \alpha_H \int_{\mathbb{R}^{d}} \int_{t}^{t+h} \int_t^{t+h} \cF G_1
(t+h-u, \cdot)(\xi) \overline{\cF G_1 (t+h-v,
\cdot)(\xi)}|u-v|^{2H-2}du \ dv \  \mu(d\xi)\\
&=& \alpha_H \int_{\mathbb{R}^{d}}\frac{\mu(d\xi)}{|\xi|^2} \int_0^h \int_0^h
\sin (u|\xi|) \sin (v|\xi|) |u-v|^{2H-2}du dv \\
&=& \alpha_H  h^{2H} \int_{\mathbb{R}^{d}}\frac{\mu(d\xi)}{|\xi|^2} \int_{0} ^{1} \int_{0} ^{1} \sin \left(u\vert \xi \vert h \right)
\sin \left( v\vert \xi \vert h \right) \vert u-v\vert ^{2H-2} dudv \\
&=& \alpha _{H}  h^{2H+2 -\beta} \int_{\mathbb{R}^{d}}\frac{\mu(d\xi)}{|\xi|^2} \int_{0} ^{1} \int_{0} ^{1}\sin (u\vert \xi \vert ) \sin (v\vert \xi \vert ) \vert u-v\vert ^{2H-2} du dv.
\end{eqnarray*}
Let us use the following notation:
\begin{equation}
\label{nt}
N_{t} (\xi)= \frac{\alpha_{H}}{\vert \xi \vert ^{2}} \int_{0} ^{t} \int_{0} ^{t} \sin (u|\xi|) \sin (v|\xi|) |u-v|^{2H-2}du dv, \hskip0.3cm t\in [0,T], \xi \in \mathbb{R} ^{d} .
\end{equation}
By Proposition 3.7 in \cite{BT} the term
$$N_{1}(\xi) = \frac{\alpha_{H}}{|\xi|^2} \int_0^1 \int_0^1
\sin (u|\xi|) \sin (v|\xi|) |u-v|^{2H-2}du dv $$ 
satisfies the inequality
\begin{eqnarray*}
N_1(\xi) & \leq & C_{H}\left(\frac{1}{1+|\xi|^2}
\right)^{H+1/2},
\end{eqnarray*}
with $C_{H}$  a positive constant not depending on $h$. Consequently the  term $E_{2}(h) $  is bounded by
\begin{equation}\label{rel2}
E_{2}(h) \leq C h^{2H+2 -\beta }   \int_{\mathbb{R} ^{d}}\left( \frac{1}{1+ \vert \xi \vert ^{2}} \right) ^{H+ \frac{1}{2}}\mu (d\xi)
\end{equation}
and this is clearly finite due to (\ref{b1}). Relations (\ref{rel3}) and (\ref{rel2}) give the first part of the conclusion.
\\
\hspace*{1.1cm} Let us analyze now the lower bound of the increments of $u(t,x)$ with respect to the variable $t$. Let $h>0, x\in [-M, M ] ^{d}$ and $t\in [t_{0}, T]$ such that
$t+h\in [t_{0}, T]$.
From the decomposition
\begin{eqnarray*}
\mathbf{E}\left| u(t+h, x) -u(t,x) \right| ^{2} &=& \Vert \left( g_{t+h, x} -g_{t,x} \right) 1_{[0,t] } \Vert _{{\cal{HP}}} ^{2} + \Vert g_{t+h, x} 1_{[t,t+h] } \Vert _{{\cal{HP}}}^{2} \\
&&+ 2 \langle \left( g_{t+h, x} -g_{t,x} \right) 1_{[0,t] }, g_{t+h, x} 1_{[t,t+h] } \rangle _{{\cal{HP}}}
\end{eqnarray*}
we immediately obtain, since the second summand in the right-hand side is positive,
\begin{eqnarray*}
\mathbf{E}\left| u(t+h, x) -u(t,x) \right| ^{2} &\geq & \Vert \left( g_{t+h, x} -g_{t,x} \right) 1_{[0,t] } \Vert _{{\cal{HP}}} ^{2} + 2 \langle \left( g_{t+h, x} -g_{t,x} \right) 1_{[0,t] }, g_{t+h, x} 1_{[t,t+h] } \rangle _{{\cal{HP}}}\\
&:=& E_{1,t}(h)+ E_{3,t}(h).
\end{eqnarray*}
We can assume, without any loss of the generality, that $t=\frac{1}{2}$. Denote $E_{1,\frac{1}{2}}(h):= E_{1}(h)$. We first prove that
\begin{equation}
\label{rel4}
E_{1} (h) \geq c h^{2H+1-\beta}-c' h^{2H+2-\beta}.
\end{equation}
for $h$ small enough. Recall that we have an exact expression for $E_{1}(h)$ (see (\ref{n11})). Actually,
\begin{eqnarray*}
E_{1}(h)&=& \int_{\mathbb{R}^{d}} \frac {d\xi } {\vert \xi \vert ^{d-\beta +2H+2}} \sin (h\vert \xi \vert ) ^{2} \Vert \cos (\cdot) 1_{(h\vert \xi \vert , h\vert \xi \vert + \vert \xi \vert )} \Vert ^{2} _{\cal{H}} \\
&=&\alpha _{H}\int_{\mathbb{R} ^{d}}
  \frac{d\xi}{\vert \xi \vert ^{d-\beta +2H+2 }}  \sin (h\vert \xi \vert ) ^{2} \int_{h\vert \xi \vert } ^{(1+h) \vert \xi \vert }  \int_{h\vert \xi \vert } ^{(1+h) \vert \xi \vert } dudv \vert u-v\vert ^{2H-2} \cos u \cos v \\
  &=& \alpha _{H}\int_{\mathbb{R} ^{d}}
  \frac{d\xi}{\vert \xi \vert ^{d-\beta +2H+2 }}  \sin (h\vert \xi \vert ) ^{2}  \int_{0} ^{\vert \xi \vert } \int_{0} ^{\vert \xi \vert } dudv \cos (u+ h\vert \xi \vert) \cos (v+ h\vert \xi \vert ) \vert u-v\vert ^{2H-2}.
\end{eqnarray*}
By the trigonomtric formula $\cos (x+y)= \cos (x) \cos (y)-\sin (x) \sin(y)$ we can write
\begin{eqnarray*}
E_{1}(h)&=& \int_{\mathbb{R}^{d}} \frac {d\xi } {\vert \xi \vert ^{d-\beta +2H+2}} \sin (h\vert \xi \vert ) ^{2}
\left[ \cos (h\vert \xi \vert ) ^{2} \int_{0} ^{\vert \xi \vert } \int_{0} ^{\vert \xi \vert } dudv \cos u \cos v \vert u-v \vert ^{2H-2}\right. \\
 &&\left. -2 \sin (h\vert \xi \vert ) \cos (h\vert \xi \vert ) \int_{0} ^{\vert \xi \vert } \int_{0} ^{\vert \xi \vert } dudv \sin u \cos v \vert u-v \vert ^{2H-2}\right. \\
 &&\left.
+ \sin (h\vert \xi \vert ) ^{2} \int_{0} ^{\vert \xi \vert } \int_{0} ^{\vert \xi \vert } dudv \sin u \sin v \vert u-v \vert ^{2H-2}\right]\\
&:=& A+B+C.
\end{eqnarray*}

We will neglect the first term since it is positive. We will bound the second one above  by $c h^{2H+2-\beta}.$ We have (we use Lemma \ref{10m-1} at the third line below)
\begin{eqnarray*}
&&\sin (h\vert \xi \vert ) \cos (h\vert \xi \vert ) \int_{0} ^{\vert \xi \vert } \int_{0} ^{\vert \xi \vert } dudv \sin u \cos v \vert u-v \vert ^{2H-2}\\
&\leq & \left| \sin (2h\vert \xi \vert )\right| \left| \int_{0} ^{\vert \xi \vert } \int_{0} ^{\vert \xi \vert } dudv \sin u \cos v \vert u-v \vert ^{2H-2}\right| \\
&= & \frac{1}{2}\left| \sin (2h\vert \xi \vert )\right| \left| \int_{0} ^{\vert \xi \vert } \int_{0} ^{\vert \xi \vert } dudv \sin (u+v)\vert u-v \vert ^{2H-2}\right| \\
&=& c \left| \sin (2h\vert \xi \vert )\right| \left| \int_{0} ^{\vert \xi \vert } v^{2H-2}(\cos (v) - \cos (2\vert \xi \vert -v)) dv \right|
\end{eqnarray*}
and thus
\begin{eqnarray*}
B&\leq & c\int_{\mathbb{R}^{d}} \frac {d\xi } {\vert \xi \vert ^{d-\beta +2H+2}} \sin (h\vert \xi \vert ) ^{2} \left| \sin (2h\vert \xi \vert )\right|
\left| \int_{0} ^{\vert \xi \vert }  v^{2H-2}(\cos (v) - \cos (2\vert \xi \vert -v) )dv \right|\\
&\leq & c \int_{\mathbb{R}^{d}} \frac {d\xi } {\vert \xi \vert ^{d-\beta +2H+2}} \sin (h\vert \xi \vert ) ^{2} \left| \sin (2h\vert \xi \vert ) \sin (\vert \xi \vert )\right|\left|    \int_{0} ^{\vert \xi \vert }\sin (v+ \vert \xi \vert )  v^{2H-2}dv\right| \\
&=&c \int_{\mathbb{R}^{d}} \frac {d\xi } {\vert \xi \vert ^{d-\beta +2H+2}} \sin (h\vert \xi \vert ) ^{2} \left| \sin (2h\vert \xi \vert ) \sin (\vert \xi \vert )\right| \\
&&\times \left| \left(  \sin (\vert \xi \vert )\int_{0} ^{\vert \xi \vert } \cos v v^{2H-2}dv +  \cos (\vert \xi \vert )\int_{0} ^{\vert \xi \vert } \sin v v^{2H-2}dv\right) \right|\\
&=&c \int_{\vert \xi \vert \leq 1} \frac {d\xi } {\vert \xi \vert ^{d-\beta +2H+2}} \sin (h\vert \xi \vert ) ^{2} \left| \sin (2h\vert \xi \vert ) \sin (\vert \xi \vert )\right| \\
&&\times \left| \left(  \sin (\vert \xi \vert )\int_{0} ^{\vert \xi \vert } \cos v v^{2H-2}dv +  \cos (\vert \xi \vert )\int_{0} ^{\vert \xi \vert } \sin v v^{2H-2}dv\right) \right|\\
&+& c \int_{\vert \xi \vert \geq 1} \frac {d\xi } {\vert \xi \vert ^{d-\beta +2H+2}} \sin (h\vert \xi \vert ) ^{2} \left| \sin (2h\vert \xi \vert ) \sin (\vert \xi \vert )\right| \\
&&\times \left| \left(  \sin (\vert \xi \vert )\int_{0} ^{\vert \xi \vert } \cos v v^{2H-2}dv +  \cos (\vert \xi \vert )\int_{0} ^{\vert \xi \vert } \sin v v^{2H-2}dv\right) \right|.
\end{eqnarray*}
The part over the set $\vert \xi \vert \leq 1$ is bounded by $c h^{3}$ by simply majorizing $\sin (x)$ by $x$ and  $\left| \left(  \sin (\vert \xi \vert )\int_{0} ^{\vert \xi \vert } \cos v v^{2H-2}dv +  \cos (\vert \xi \vert )\int_{0} ^{\vert \xi \vert } \sin v v^{2H-2}dv\right) \right|$ by a constant.
Concerning the part over the region $\vert \xi \vert \geq 1$ we bound again $\left| \left(  \sin (\vert \xi \vert )\int_{0} ^{\vert \xi \vert } \cos v v^{2H-2}dv +  \cos (\vert \xi \vert )\int_{0} ^{\vert \xi \vert } \sin v v^{2H-2}dv\right) \right|$ by a constant and we use the change of variables $\tilde{\xi} = h\xi$. This part will by bounded by
\begin{eqnarray*}
&&h^{2H+2-\beta}\int_{\vert \xi \vert \geq h} \frac {d\xi } {\vert \xi \vert ^{d-\beta +2H+2}} \sin (\vert \xi \vert ) ^{2} \left| \sin (2\vert \xi \vert ) \right| \\
&&\leq h^{2H+2-\beta}\int_{\mathbb{R} ^{d}} \frac {d\xi } {\vert \xi \vert ^{d-\beta +2H+2}} \sin (\vert \xi \vert ) ^{2} \left| \sin (2\vert \xi \vert ) \right|\\
&&\leq ch^{2H+2-\beta }
\end{eqnarray*}
since the last integral is convergent at infinity by bounded sinus by one and at zero by bounding $\sin (x)$ by $x$ and using the assumption $\beta >2H-1$. Therefore
\begin{equation}
\label{rel5}
B\leq ch^{2H+2-\beta}.
\end{equation}

We bound now the summand $C$ below. In this summand the ${\cal{H}}$ norm of the sinus function appear and this has been analyzed in \cite{BT}. We have, after the change of variables $\tilde{u}= \frac{u}{\vert \xi \vert }, \tilde{v}=\frac{v}{\vert \xi \vert }$,
\begin{eqnarray*}
C&=& \int_{\mathbb{R}^{d}} \frac {d\xi } {\vert \xi \vert ^{d-\beta +2}} \sin (h\vert \xi \vert ) ^{4}\int _{0} ^{1} \int_{0} ^{1} \sin (\vert \xi \vert) \sin (v\vert \xi \vert ) \vert u-v\vert ^{2H-2} dudv \\
&\geq & \int_{\vert \xi \vert \geq 1} \frac {d\xi } {\vert \xi \vert ^{d-\beta +2}} \sin (h\vert \xi \vert ) ^{4}\int _{0} ^{1} \int_{0} ^{1} \sin (\vert \xi \vert) \sin (v\vert \xi \vert ) \vert u-v\vert ^{2H-2} dudv.
\end{eqnarray*}
We will use Proposition 3.8 in \cite{BT} (more precisely, we will use the inequality two lines before the formula (34) in that paper with $k=0$; we notice that the term $\sin (h\vert \xi \vert ) ^{4}$ does not appear in this proof but by analyzing the step of the proof we can see that this term can be added without problems). We will have that, for $h$ small,
\begin{eqnarray}
C&\geq &  \int_{\vert \xi \vert \geq 1} \frac {d\xi } {\vert \xi \vert ^{d-\beta }} \sin (h\vert \xi \vert ) ^{4}\frac{1}{\vert \xi \vert ^{2} }\int _{0} ^{1} \int_{0} ^{1} \sin (\vert \xi \vert) \sin (v\vert \xi \vert ) \vert u-v\vert ^{2H-2} dudv \nonumber \\
&\geq & \int_{\vert \xi \vert \geq 1} \frac {d\xi } {\vert \xi \vert ^{d-\beta }} \sin (h\vert \xi \vert ) ^{4} \frac{1}{\vert \xi \vert ^{2H+1}}\nonumber \\
&=& h^{2H+1-\beta} \int _{\vert \xi \vert \geq h} \frac {d\xi } {\vert \xi \vert ^{d-\beta  +2H+1 }} \sin (h\vert \xi \vert ) ^{4}\nonumber  \\
&\geq & h^{2H+1-\beta} \int _{\vert \xi \vert \geq 1} \frac {d\xi } {\vert \xi \vert ^{d-\beta  +2H+1 }} \sin (h\vert \xi \vert ) ^{4}\nonumber \\
&=& ch^{2H+1-\beta}.\label{rel6}
\end{eqnarray}

  Relations (\ref{rel5}) and (\ref{rel6}) imply (\ref{rel4}). Now, from relation (\ref{rel4}), for every $t_{0} \leq s<t <T$ with $s,t$ close enough
  \begin{eqnarray*}
  E_{1} (t-s) &\geq & c (t-s)^{2H+1-\beta}-c' (t-s)^{2H+2-\beta} \geq \frac{c}{2} (t-s)^{2H+1-\beta}
\end{eqnarray*}
if $\vert t-s\vert \leq \frac{c}{2c'}.$ To extend the above inequality to arbitrary values of $\vert t-s\vert$, we proceed  as in \cite{DaSa1}, proof of Proposition 4.1. Notice that the function $g(t,s, x,y):= \mathbf{E}\left| u(t,x)-u(s,x) \right| ^{2}$ is positive and continuous with respect to all its arguments and therefore it is bounded below on the set $\{ (t,s,x,y) \in [t_{0}, T] ^{2} \times [-M, M] ^{2d}; \vert t-s\vert \geq \varepsilon \} $ by a constant depending on $\varepsilon >0$. Hence for $\vert t-s\vert \geq \frac{c}{2c'}$ it also holds that
\begin{equation*}
E_{1} (t-s)\geq c_{1} \vert t-s\vert ^{2H+1-\beta}.
\end{equation*}

 On the other side,  from (\ref{rel1}) and (\ref{rel2}) and Cauchy-Schwarz inequality, we obtain
\begin{eqnarray*}
E_{3,t}(h)&=& \langle \left( g_{t+h, x} -g_{t,x} \right) 1_{[0,t] }, g_{t+h, x} 1_{[t,t+h] } \rangle _{{\cal{HP}}}\\
&&\leq \Vert \left( g_{t+h, x} -g_{t,x} \right) 1_{[0,t] } \Vert _{{\cal{HP}}} \Vert g_{t+h, x} 1_{[t,t+h] } \Vert _{{\cal{HP}}}\\
&&\leq c h^{\frac{2H+1-\beta}{2}  + \frac{2H+2-\beta}{2}}.
\end{eqnarray*}
Consequently,
\begin{eqnarray*}
\mathbf{E}\left| u(t+h, x) -u(t,x) \right| ^{2} &\geq & C h^{2H+1-\beta} -C'  h^{\frac{2H+1-\beta}{2} + \frac{2H+2-\beta}{2}}
\end{eqnarray*}
and this implies that for every $s,t\in [t_{0}, T]$ and $x\in [-M, M ] ^{d}$
\begin{eqnarray*}
\mathbf{E}\left| u(t,x)-u(s,x) \right| ^{2}&\geq & \frac{C}{2}\vert t-s\vert ^{2H + 1 -\beta} \quad \mbox{ if } \quad \vert t-s\vert \leq \left( \frac{C}{2C'} \right) ^{ \frac{1}{2 }}.
\end{eqnarray*}
Similarly as above, the previous inequality can be extended to arbitrary values of $s,t\in [t_{0} , T]$.

\qed

\vskip0.2cm
Proposition (\ref{ptime}) implies the following H\"older property for the solution to (\ref{wave}).
 \begin{corollary}
 Assume (\ref{b2}). Then for every $x \in \mathbb{R} ^{d}$  the application
$$t\to u(t,x)$$
is almost surely  H\"older continuous of order $\delta \in \left( 0, \frac{2H+1-\beta}{2}\right). $
\end{corollary}
{\bf Proof: } This is consequence of the relations (\ref{rel1}) and (\ref{rel2}) in the proof of Proposition \ref{ptime} and of the fact that $u$ is Gaussian. \qed

\vskip0.2cm

Let us make some comments on the result in Proposition \ref{ptime}.
\begin{remark}
\begin{description}
\item
\item{$\bullet$ } Following the proof of Theorem 5.1 in \cite{DaSa2} we can show that the mapping $t\to u(t,x)$ is not H\"older continuous of order   $\frac{2H+1-\beta}{2}$.
 \item{$\bullet$ } When $H=\frac{1}{2}$  the solution of the wave equation with fractional noise in time has the same regularity in time as the solution of the wave equation with white noise in time (see \cite{DaSa1}, \cite{DaSa2}).
     \end{description}
  \end{remark}

\subsection{Space regularity}
\hspace{1.1cm} Let us discuss the behavior of the solution $u$ to the equation (\ref{wave}) with respect to the spatial variable. We have


\begin{prop}\label{pspace}
 Assume (\ref{b2}), fix $M>0$ and $t\in [t_{0}, T]$. Then there exist positive constants $c_{3}, c_{4}$ such that for any $x,y \in [-M, M] ^{d} $
\begin{equation*}
c_{3} \vert x-y\vert ^{2H+1-\beta} \leq \mathbf{E} \left| u(t,x) -u(t,y) \right| ^{2} \leq c_{4} \vert x-y\vert ^{2H+1-\beta }.
\end{equation*}

\end{prop}
{\bf Proof: } Let $z\in \mathbb{R}^d $. We compute
\begin{eqnarray*}
&&\mathbf{E}|u(t,x+z)-u(t,x)|^2 =\|g_{t,x+z}-g_{t,x}\|_{\cH \cP}^2 \\
& & = \alpha_H \int_{\mathbb{R}^{d}} \int_{0}^{t} \int_{0}^{t} \cF
(g_{t,x+z}-g_{t,x}) (u, \cdot)(\xi) \overline{\cF (g_{t,x+z}-g_{t,x})(v,
\cdot)(\xi)}|u-v|^{2H-2}du \ dv \  \mu(d\xi) \\
& & = \alpha_H \int_{0}^{t} \int_{0}^{t} |u-v|^{2H-2}du \ dv  \int_{\mathbb{R}^{d}} |e^{-i \xi \cdot (x+z)}-e^{-i \xi \cdot
x}|^2 \cF G_1(u, \cdot)(\xi) \overline{\cF
G_1(v, \cdot)(\xi)}\  \mu(d\xi) \\
& & = \alpha_H \int_{0}^{t} \int_{0}^{t} |u-v|^{2H-2}du \ dv  \int_{\mathbb{R}^{d}} |e^{-i \xi \cdot z}-1|^2 \frac{\sin (u|\xi|)}{|\xi|} \cdot \frac{\sin (v|\xi|)}{|\xi|} \mu(d\xi)\\
& & =: E_{1,x}(z) + E_{2,x}(z),
\end{eqnarray*}
where $E_{1,x}(z)$ and $E_{2,x}(z)$ are the integrals over the regions $ \vert \xi \vert < 1$ and $\vert \xi \vert \geq 1$ respectively.\\
For the first expresion is easy to see that, using the inequality $\vert 1-e^{-i\xi z} \vert ^{2} \leq \vert \xi \vert ^{2} \vert z \vert  ^{2}$, we get the bound
\begin{equation*}
	E_{1,x}(z) \leq C \vert z \vert  ^{2} \int _{\vert \xi \vert \leq 1} \mu (d\xi).
\end{equation*}
Developing the second expresion we get
\begin{eqnarray*}
E_{2,x}(z) &=& \alpha_H \int_{0}^{t} \int_{0}^{t} |u-v|^{2H-2}du \ dv  \int_{\vert \xi \vert \geq 1} |e^{-i \xi \cdot z}-1|^2 \frac{\sin (u|\xi|)}{|\xi|} \cdot \frac{\sin (v|\xi|)}{|\xi|} \mu(d\xi)\\
&=& 2\alpha_H \int_{0}^{t} \int_{0}^{t} |u-v|^{2H-2}du \ dv  \int_{\vert \xi \vert \geq 1}\frac{d\xi}{\vert \xi \vert ^{d-\beta}} (1-\cos (z\cdot \xi))\frac{\sin (u|\xi|)}{|\xi|} \cdot \frac{\sin (v|\xi|)}{|\xi|},
\end{eqnarray*}
where $z\cdot \xi $ means the scalar product in $\mathbb{R} ^{d}$. Again from Proposition 3.7 in \cite{BT} we have that
 $$N_t(\xi) \leq c_{t,H} \left(\frac{1}{1+|\xi|^2}\right)^{H+1/2} $$ for any $t>0$, $ |\xi| \geq 1 $, where $N_{t}(\xi)$ is given by (\ref{nt}). Hence, denoting by $e= \frac{z}{\vert z \vert }$

 \begin{eqnarray*}
E_{2,x}(z) &\leq & C \int_{\mathbb{R}^{d}}\frac{d\xi}{\vert \xi \vert ^{d-\beta}} (1-\cos (z\cdot \xi))\left( \frac{1}{1+ \vert \xi \vert ^{2}}\right) ^{H+\frac{1}{2}}\\
&=& C z^{2H+1-\beta } \int_{\mathbb{R}^{d}} \frac{dw} {\vert w\vert ^{d-\beta }}(1-\cos (w \cdot e)) \left( \frac{1}{\vert  w\vert ^{2} + \vert z\vert ^{2} } \right) ^{H+\frac{1}{2}}\\
&\leq & C \vert z\vert ^{2H+1-\beta},
 \end{eqnarray*}
 where we used the change of variables $w= \xi \vert z\vert  $. This proves the upper bound.

Let us prove the sharpness of this bound (i.e. the lower bound). We can assume, without losing the generality, that $t=1$. We note that
\begin{eqnarray*}
&&\mathbf{E} \left| u(1, x+z)- u(1,x) \right| ^{2} \geq  F_{2} (z) \\
&:=& 2 \alpha _{H} \int_{0} ^{1} \int_{0} ^{1} dudv \vert u-v\vert ^{2H-2} \int _{\vert \xi \vert \geq 1} \frac{d\xi } {\vert \xi \vert ^{d-\beta } } (1-\cos ( \xi \cdot z)) \frac{\sin (u|\xi|)}{|\xi|} \cdot \frac{\sin (v|\xi|)}{|\xi|}.
\end{eqnarray*}

Condition (\ref{mu-tempered-eq}) implies that
\begin{equation*}
\int_{\vert \xi \vert \geq 1} \frac{\mu (d\xi) }{\vert \xi \vert ^{3} } \leq  \int_{\vert \xi \vert \geq 1} \frac{\mu (d\xi) }{\vert \xi \vert ^{2H+1} }<\infty.
\end{equation*}

We apply Proposition 3.8 in \cite{BT} (more precisely, the inequality two lines before (34) in \cite{BT} with $k=0$)  and we get (note that the result in \cite{BT} is stated without the factor $(1-\cos ( \xi \cdot z))$ but by analyzing the steps of the proof we can see that this factor may be added without problems)
\begin{equation*}
F_{2} (z)\geq C \int_{\vert \xi \vert \geq 1} \frac{d \xi } {\vert \xi \vert ^{d-\beta +2H+1} } (1-\cos ( \xi \cdot z) ).
\end{equation*}
and by the change of variables $\xi \vert z \vert =w$ in the integral $d\xi$
\begin{eqnarray}
\nonumber
F_{2} (z) &\geq & C z^{2H+1-\beta } \int_{\vert w \vert \geq  \vert z \vert}\frac{dw}{\vert w\vert ^{d-\beta +2H +1}}(1-\cos (w \cdot e).
\label{n1}
\end{eqnarray}
As in the proof of Theorem 5.1 in \cite{DaSa2}, we obtain that the integral $\int_{\vert w \vert \geq  \vert z \vert}\frac{dw}{\vert w\vert ^{d-\beta +2H +1}}(1-\cos (w\cdot e))$ is bounded below by a constant. (Notice that $\beta >2H-1$, implies that the first integral above is convergent when $z$ is zero, because $1-\cos (x) \thickapprox x^{2} $ around zero).   Thus, it is immediate that
\begin{equation*}
\mathbf{E} \left| u(1, x+z)- u(1,x) \right| ^{2} \geq  C z^{2H+1-\beta} -C' z^{2H+2-\beta }
\end{equation*}
and this implies
\begin{equation*}
\mathbf{E} \left| u(1, x+z)- u(1,x) \right| ^{2} \geq  \frac{1}{2}C z^{2H+1-\beta}
\end{equation*}
for $\vert z\vert \leq \frac{C}{2C'}$. It is a routine argument to extend the above inequality to arbitrary values of $\vert z\vert$ (see e.g. \cite{DaSa1}, page 22, see also the proof of Proposition \ref{ptime} before).

\qed

\vskip0.2cm

We have the following result concerning the H\"older continuity in space. We mention that it is a little bit more  than  an extension of Proposition \ref{pspace}.

\begin{prop} \label{p2}
Assume (\ref{b2}). Then for any $t\in [t_{0}, T]$ the application
$$x\to u(t,x)$$
is almost surely H\"older continuous of order $\delta \in \left(0, \left( \frac{2H+1-\beta}{2} \right) \wedge 1 \right)$.
\end{prop}
{\bf Proof: }
We claim that
\begin{equation}\label{rel7}
\mathbf{E} \left| u(t,x) -u(t,y) \right| ^{2} \leq c \vert x-y\vert ^{(2H+1-\beta) \wedge 2 }
\end{equation}
for $|x-y|$ small  enough.  From Proposition \ref{pspace}, (\ref{rel7}) is true when $\beta >2H-1$.
When $\beta \in (0, 2H-1)$ then it suffices to regards the part of the quantity $\mathbf{E}\left| u(t,x+z)-u(t,x)\right| ^{2} $ over the region $\vert \xi \vert \leq 1$  (the part over the region $\vert \xi \vert >1$ is, as in the proof of Proposition \ref{pspace}, bounded by $cz^{2H+1-\beta}$ so by $cz^{2}$ for $z$ small). It is immediate to see that, using the inequality $\vert 1-e^{-i\xi z} \vert ^{2} \leq \vert \xi \vert ^{2} \vert z\vert ^{2}$ the considered part is less than $C \vert z\vert ^{2} \int _{\vert \xi \vert \leq 1} \mu (d\xi). $ This concludes the proof of (\ref{rel7}).

The conclusion  is a consequence of  Proposition \ref{pspace}, the Gaussianity of $u$  and the Kolmogorov continuity theorem.

\qed
\vskip0.2cm

\begin{remark}

\begin{description}
\item
\item {$\bullet$ } When $H=\frac{1}{2}$  the solution of the wave equation with fractional noise in time has the same regularity in space as the solution of the wave equation with white noise in time (see \cite{DaSa1}, \cite{DaSa2}).

\item{$\bullet$ }We distinguish in Proposition \ref{p2} two cases: if $\beta \in (0, 2H-1)$ then the solution to (\ref{wave})  has spatial H\"older continuity of order 1  (so, it is Lipschitz continuous in the space variable) while if $\beta \in (2H-1, d\wedge (2H+1))$ the H\"older exponent is $\delta \in (0, \frac{2H+1-\beta}{2})<1$.

\item{$\bullet$ }There is another way to see why the cases $\beta \in (0,2H-1]$ and $\beta \in  (2H-1,  d\wedge (2H+1))$ need to be separated. Denote by
\begin{eqnarray*}
g_{t}(z) &:=& \mathbf{E} \left| u(t, x+z) -u(t,x) \right| ^{2} \\
&=&  2\alpha_H \int_{0}^{t} \int_{0}^{t} |u-v|^{2H-2}du \ dv  \int_{\mathbb{R}^{d}}\frac{d\xi}{\vert \xi \vert ^{d-\beta}} (1-\cos (z\cdot \xi))\frac{\sin (u|\xi|)}{|\xi|} \cdot \frac{\sin (v|\xi|)}{|\xi|}
\end{eqnarray*}
and let us study the behavior of $g_{t}$ around $z=0$. Let us also assume that $d=1$. Notice first that $g_{t}(0)=0$ and

\begin{equation*}
g'_{t}(z)= 2\alpha_H \int_{0}^{t} \int_{0}^{t} |u-v|^{2H-2}du \ dv  \int_{\mathbb{R}^{d}}\frac{d\xi}{\vert \xi \vert ^{d-\beta -1}} \sin (z\cdot \xi )\frac{\sin (u|\xi|)}{|\xi|} \cdot \frac{\sin (v|\xi|)}{|\xi|}
\end{equation*}

and thus $g'_{t}(0)=0$ provided that $\beta <2H$. Moreover
\begin{equation*}
g''_{t}(z)=  2\alpha_H \int_{0}^{t} \int_{0}^{t} |u-v|^{2H-2}du  dv  \int_{\mathbb{R}^{d}}\frac{d\xi}{\vert \xi \vert ^{d-\beta-2}} \cos (z\cdot \xi)\frac{\sin (u|\xi|)}{|\xi|} \cdot \frac{\sin (v|\xi|)}{|\xi|}
\end{equation*}

and
\begin{eqnarray*}
g''_{t}(0)&= &2\alpha_H \int_{0}^{t} \int_{0}^{t} |u-v|^{2H-2}du  dv  \int_{\mathbb{R}^{d}}\frac{d\xi}{\vert \xi \vert ^{d-\beta-2}}  \frac{\sin (u|\xi|)}{|\xi|} \cdot \frac{\sin (v|\xi|)}{|\xi|}\\
&\leq & C_{t} 2\alpha_H \int_{\mathbb{R}^{d}}\frac{d\xi}{\vert \xi \vert ^{d-\beta-2}} \left( \frac{1}{1+ \vert \xi \vert ^{2}} \right) ^{H+\frac{1}{2}}
\end{eqnarray*}

which is a finite constant for $\beta <2H-1$. Therefore $g_{t}(z)$  behaves as $Cz^{2}$ for $z$ close to zero.

\end{description}

\end{remark}

\subsection{Joint regularity}

\hspace{1.1cm} Let us denote by $\Delta $ the following metric on $[0,T] \times \mathbb{R} ^{d}$
\begin{equation}
\label{metric}
\Delta \left( (t,x); (s,y) \right)= \vert t-s\vert ^{2H+1-\beta } + \vert x-y\vert^{ 2H+1-\beta }.
\end{equation}
From Propositions \ref{ptime} and \ref{pspace}, we obtain the following result:

\begin{theorem}
\label{increm}
Fix $M>0$ and assume (\ref{b2}). For every $t,s \in [t_{0}, T]$ and $x,y\in [-M, M] ^{d}$ there exist positive constants $C_{1}, C_{2}$ such that
\begin{equation*}
C_{1}\Delta \left( (t,x); (s,y) \right) \leq \mathbf{E} \left| u(t,x) -u(s,y) \right| ^{2} \leq C_{2} \Delta \left( (t,x); (s,y) \right).
\end{equation*}
\end{theorem}
{\bf Proof: } The upper bound can be easily obtained by using the upper bound in Propositions \ref{ptime} and \ref{pspace} since
\begin{eqnarray*}
\mathbf{E} \left| u(t,x) -u(s,y) \right| ^{2}& \leq & 2 \mathbf{E} \left| u(t,x) -u(s,x) \right| ^{2}+ 2\mathbf{E} \left| u(s,x) -u(s,y) \right| ^{2}\\
&\leq & C_{2} \left( \vert t-s\vert  ^{2H+1-\beta} + \vert x-y\vert ^{2H+1-\beta }\right).
\end{eqnarray*}

Concerning the lower bound, it suffices follow the lines of the proof of Lemma 2.1 in \cite{EF} (see also Steps 3 and 4 in the proof of Proposition 4.1 in \cite{DaSa1}).  We will briefly explain the main lines of the proof. The demostration needs to be divided upon three cases: $ \vert t-s\vert ^{2H+1-\beta } \leq \frac{c_{3}}{4c_{2}} \vert x-y\vert ^{2H+1-\beta}$,\ \  $ \vert t-s\vert ^{2H+1-\beta } \geq \frac{4c_{4}}{c_{1}} \vert x-y\vert ^{2H+1-\beta}$ and $\frac{4c_{4}}{c_{1}} \vert x-y\vert ^{2H+1-\beta} \geq \vert t-s\vert ^{2H+1-\beta }  \geq \frac{c_{3}}{4c_{2}} \vert x-y\vert ^{2H+1-\beta}$ with the constants $c_{1}, c_{2}, c_{3}, c_{4}$ appearing in the statements of Propositions \ref{ptime} and \ref{pspace}. The first case can be handled as follows
\begin{eqnarray*}
\mathbf{E} \left| u(t,x) -u(s,y) \right| ^{2}& \geq & \frac{1}{2} \mathbf{E} \left| u(t,x) -u(t,y) \right| ^{2}- \mathbf{E} \left| u(t,y) -u(s,y) \right| ^{2}\\
& \geq & \frac{1}{2} c_{3}\vert x-y\vert ^{2H+1-\beta} - c_{2} \vert t-s\vert ^{2H+1-\beta} \\
&\geq & \frac{1}{2} c_{3}\vert x-y\vert ^{2H+1-\beta} -\frac{1}{4} c_{3} \vert x-y\vert ^{2H+1-\beta }\\
&= & \frac{1}{4} c_{3} \vert x-y\vert ^{2H+1-\beta}\\
&\geq & \frac{c_{3}}{8} \vert x-y\vert ^{2H+1-\beta}+ \frac{c_{3}}{8}  \frac{4c_{2}}{c_{3} } \vert t-s\vert ^{2H+1-\beta} \\
&\geq &C_{1} \Delta \left( (t,x); (s,y) \right).
\end{eqnarray*}

The other cases follows similarly from Lemma 3.1  in \cite{EF}, by replacing their exponents with our exponents.

\qed

\vskip0.2cm

\begin{remark}
The result of Theorem \ref{increm} can be stated also in the following form: Fix $M>0$ and assume (\ref{b2}). For every $t,s \in [t_{0}, T]$ and $x,y\in [-M, M] ^{d}$ with $(t,x)$ close enough to $(s,y)$, there exist positive constants $C_{1}, C_{2}$ such that
\begin{equation*}
C_{1}\left( \vert t-s\vert + \vert x-y\vert \right) ^{2H+1-\beta }  \leq \mathbf{E} \left| u(t,x) -u(s,y) \right| ^{2} \leq C_{2} \left( \vert t-s\vert + \vert x-y\vert \right) ^{2H+1-\beta}.
\end{equation*}
\end{remark}

\vspace*{0.2cm}

\section{Hitting times}

\hspace{1.1cm} Let us discuss the upper and lower bounds for the hitting probabilities of the solution $u$ to equation (\ref{wave}). These bounds will be given   in terms of the Newtonian capacity and the Hausdorff measure of the hit set (see Section 2 for the definition). Let us recall the notation: if $V=(V(x), x\in \mathbb{R} ^{m})$ is a $\mathbb{R} ^{k}$ valued stochastic process then $V(S)$ denote the range of the Borel set $S$ under the random mapping $x\to V(x)$.

Our result is based on the following criteria for the hitting probabilities proven in \cite{BLX}, Theorem 2.1.

\begin{theorem}
\label{blx}
Let $X= X(t), t\in \mathbb{R} ^{N}$ be a $\mathbb{R} ^{k}$-valued centered Gaussian process  and fix $I\subset \mathbb{R} ^{N}$. Assume that there exist positive constants $a_{1}, a_{2} , a_{3} , a_{4}$ such that
\begin{description}
\item{i. } For every $t\in I$, $\mathbf{E} \left[ X(t)^{2}\right] \geq a_{1}>0  $.

\item{ii. } There exists $\alpha_{1},..., \alpha _{N} \in (0,1)$ such that for every $t=(t_{1},...,t_{N}), s=(s_{1},..., s_{N}) \in I $ it holds that
    \begin{equation*}
    a_{2}\sum_{j=1} ^{N} \vert t_{j}-s_{j} \vert ^{2\alpha _{j}} \leq \mathbf{E} \left| X(t)-X(s) \right| ^{2} \leq a_{3}\sum_{j=1} ^{N} \vert t_{j}-s_{j} \vert ^{2\alpha _{j}}.
    \end{equation*}
\item{iii. }For every $t=(t_{1},...,t_{N}), s=(s_{1},..., s_{N}) \in I $
\begin{equation*}
Var (X(t)| X(s) ) \geq a_{4}\sum_{j=1} ^{N} \vert t_{j}-s_{j} \vert ^{2\alpha _{j}}.
\end{equation*}
\end{description}
Then  there exist positive constants $a_{5}, a_{6}$ such that  for every Borel set $A$ in $\mathbb{R} ^{k}$
\begin{equation*}
a_{5} \text{Cap} _{k-Q} (A) \leq P(X(I) \cap A \not= \emptyset)\leq a_{6} {\cal{H}}_{k-Q}(A)
\end{equation*}
where $Q= \sum_{j=1}^{N} \frac{1}{\alpha _{j}}.$
\end{theorem}

Next, we will show that the solution to (\ref{wave}) satisfies the assumptions of the  previous result. This will be done via several lemmas.

\begin{lemma}
Assume (\ref{b2}) and let $u$ be the solution to (\ref{wave}). Then for every $t\in [t_{0}, T]$ and $x\in \mathbb{R} ^{d}$
\begin{equation*}
\mathbf{E}u(t,x)^{2} \geq C.
\end{equation*}
{\bf Proof: } Let $\sigma ^{2}_{t,x} $ be the variance of $u(t,x)$.  We need to give a lower bound for this variance. Assume for simplicity $t=1$.
Then
\begin{eqnarray*}
\sigma ^{2} _{1,x}&=& \mathbf{E} \left| u(1,x) \right| ^{2} \\
&=&  \alpha _{H} \int_{0} ^{1} \int_{0} ^{1}dudv \vert u-v\vert ^{2H-2} \int _{\mathbb{R}^{d}} \frac{d\xi }{\vert \xi \vert ^{d-\beta +2} } \sin (u\vert \xi \vert ) \sin (v\vert \xi \vert ) \\
&\geq &\alpha _{H} \int_{0} ^{1} \int_{0} ^{1} dudv\vert u-v\vert ^{2H-2} \int _{\vert \xi \vert \leq 1} \frac{d\xi }{\vert \xi \vert ^{d-\beta +2} } \sin (u\vert \xi \vert ) \sin (v\vert \xi \vert ) \\
&\geq & \alpha _{H} \sin ^{2} 1 \int_{0}^{1}  \int_{0} ^{1} dudv\vert u-v\vert ^{2H-2}uv = C>0
\end{eqnarray*}
where we used the bound $\sin x \geq x\sin 1 $ for every $x\in [0,1]$.  The general case $t\in [t_{0}, T] $ follows in the same way by doing the change of variables $\tilde{u}= \frac{u}{t}, \tilde{v}=\frac{v}{t}$ and then working on the domain $D=\{ \xi \in \mathbb{R} ^{d}, \vert \xi \vert \leq \frac{1}{ut} \}$.\qed
\end{lemma}

Now, we bound the conditional variance (condition iii. in Theorem \ref{blx}).
\begin{lemma}
Assume (\ref{b2}) and fix $t_{0},M>0$. Then for every $s,t\in [t_{0}, T]$ and $x,y\in [-M, M ] ^{d}$
\begin{equation*}
Var (u(t,x)| u(s,y) )\geq  C \Delta ((t,x); (s,y))
\end{equation*}
where $\Delta $ is the metric given by (\ref{metric}).
\end{lemma}
{\bf Proof: } We will use the following formula: if $(U,V)$ is a centered Gaussian vector, then
\begin{equation}
\label{vc}
Var(U,V) = \frac{ (\rho^{2}_{U,V} -(\sigma _{U}-\sigma _{V} )^{2}) ((\sigma _{U} + \sigma _{V}) ^{2} -\rho^{2} _{U,V} )}{4 \sigma ^{2}_{V}}
\end{equation}
where $\rho^{2}_{U,V}= \mathbf{E} (U-V) ^{2}, \sigma ^{2}_{U}= \mathbf{E}U^{2}, \sigma ^{2}_{V} = \mathbf{E}V^{2}$.
Denote by
\begin{equation*}
\rho^{2} _{t,x,s,y} = \mathbf{E}\left|  u(t,x) -u(s,y) \right| ^{2}, \hskip0.2cm \sigma ^{2}(t,x)= \mathbf{E}u(t,x)^{2}, \hskip0.2cm \sigma ^{2} _{s,y}= \mathbf{E} u(s,y) ^{2}.
\end{equation*}
It suffices to show that
\begin{equation*}
(\rho^{2} _{t,x,s,y}-(\sigma _{t,x} -\sigma _{s,y}) ^{2}) ( (\sigma _{t,x}+ \sigma _{s,y}) ^{2}-\rho^{2} _{t,x,s,y} )\geq c \Delta ((t,x); (s,y))
\end{equation*}
for every $s,t\in [t_{0}, T]$ and $x,y\in [-M, M]^{d}$. By Theorem \ref{increm} the second factor in the left-hand side above is bounded below by a constant. So it remains to check that
$$(\rho^{2} _{t,x,s,y}-(\sigma _{t,x} -\sigma _{s,y}) ^{2}) \geq c \Delta ((t,x); (s,y))$$
but this has been done in \cite{EF}, proof of Proposition 3.2. (see also \cite{DKN1}, proof of Lemma 4.3).

\qed

\vskip0.2cm

\begin{remark}
Using the  previous result we can give a bound on the joint density $p_{t,x,s,y}$ of the vector $(u(t,x), u(s,y))$. Actually, one can show that for every $t\in [t_{0}, T]$  and $x,y \in [-M,M]^{d}$ we have the inequality
\begin{equation*}
p_{t,s,x,y}(z_{1}, z_{2}) \leq C_{1} \Delta \left( (t,x); (s,y) \right) ^{-\frac{1}{2}} \exp \left( -\frac{C_{2} \vert z_{1}-z_{2} \vert ^{2}}{ \Delta \left( (t,x); (s,y) \right) }\right)
\end{equation*}
for every $z_{1}, z_{2} \in [-N,N]^{k}$, where $\Delta$ is the metric defined by (\ref{metric}). It suffices to follow the lines of Proposition 3.2 in \cite{EF}.
\end{remark}

We can state now the main result of this section.

\begin{theorem}
Let us consider $I,J$ non-trivial compacts sets in $[t_{0}, T]$ and $[-M,M]^{d}$ respectively. Fix $N>0$ and let $u$ be the solution to the system (\ref{system}).  Then for every Borel set $A$ contained in $[-N, N]^{k}$ it holds that
\begin{equation*}
C^{-1} \text{Cap}_{k-\gamma }\leq P \left( u(I\times J) \cap A \not= \emptyset\right) \leq C {\cal{H}} _{k-\gamma } (A)
\end{equation*}
with
$$\gamma = k-\frac{2(d+1) }{2H+1-\beta}.$$

\end{theorem}
{\bf Proof: } The proof is a consequence of Theorem \ref{blx} and of the preceding two lemmas.

\qed

\begin{remark}
\begin{description}
\item
\item{$\bullet$ }
Of course, for $H$ close to $\frac{1}{2}$, our result recovers the findings in \cite{DaSa1}.
\item{$\bullet$ } it is also possible to give some results concerning the probability that, for fixed $t,x$, the  sets $u(\{t\}\times J) $ and $u(I\times \{x \} )$ (as before $I,J$ non-trivial compacts sets in $[t_{0}, T]$ and in $[-M,M]^{d}$ respectively) to hit a given Borel set $A$ contained in $[-N, N]^{k}$. Actually, by routine arguments we will have
    \begin{equation*}
   C^{-1} \text{Cap}  _{k-\frac{2d}{2H+1-\beta}}(A)\leq P\left( u(\{t\}\times J) \cap A \not= \emptyset\right) \leq C {\cal{H}} _{k-\frac{2d}{2H+1-\beta}}(A)
    \end{equation*}
    and
    \begin{equation*}
     C^{-1} \text{Cap}  _{k-\frac{2}{2H+1-\beta}}(A)\leq P\left( u(I\times \{x \} ) \cap A \not= \emptyset\right) \leq C {\cal{H}} _{k-\frac{2}{2H+1-\beta}}(A).
    \end{equation*}

\end{description}
\end{remark}

{\bf Acknowledgement: } The second author is grateful to Professor Robert Dalang, Davar Khosnevisan and Marta Sanz-Sol\'e for useful discussions.

\end{document}